\documentstyle{amsppt}
\magnification1200
\NoRunningHeads
\def\>>>{\rightarrow}

\topmatter
\title
Diffeomorphisms and almost complex structures on tori
\endtitle

\author Bogus\l\/aw Hajduk$^1$ and Aleksy Tralle{}\footnote{Both authors
were partially supported by grant 2P03A 036 24 of the Polish
Committee for the Scientific Research (KBN)\ \ \ \ \ \ \ \ \ \ \ \
\ }
\endauthor

\abstract
We prove that there exist diffeomorphisms of tori, supported in a disc,
 which are not isotopic to symplectomorphisms with respect to the
 standard symplectic structure.  This yields a partial negative answer to a
question of Benson and Gordon about the existence of symplectic structures
on tori with exotic differential structure.

\endabstract
\endtopmatter
\document
\head 1. Introduction
\endhead
The initial motivation for this work was the following remark in
the paper of Benson and Gordon \cite{2}: {\it it is not known, if
there exist K\"ahler structures on exotic tori}. In fact, it is
{\it even not known if there are symplectic structures on a torus
with an exotic differential structure.} Motivated by this, we
study the possibility of building symplectic structures on  exotic
tori using the classical construction of Thurston. This
construction gives  a symplectic form on a compact manifold $M$
fibred over a symplectic  manifold with symplectic fiber provided
the fibration is symplectic and the following condition is
satisfied: there exists a class in $H^2(M,\Bbb R)$ which restricts
to the cohomology class of the symplectic form of the fibre (see
\cite{10}, 6.3).

Some exotic tori have the structure of a fiber bundle with the base
and the fiber being  even dimensional tori with standard differential structures.
This is obtained as follows. In the sequel we denote
by $\Bbb T^{k}$ the standard $k$-torus and we write ${\Cal T}^{k}$
if we consider a torus with an exotic structure.
Let $f: \Bbb T^{2n}\to\Bbb T^{2n}$ be a diffeomorphism
supported in a disc, i.e.,
 equal to the identity outside
an embedded disc $D^{2n}$. This diffeomorphism corresponds to a
diffeomorphism $\hat f$ of the $2n$-sphere $S^{2n}$ and hence
gives  a $(2n+1)$-dimensional homotopy sphere
$\Sigma_f=D^{2n+1}\cup_{\hat f}D^{2n+1}.$ Simply consider the disc
were $f$ is supported as embedded in the sphere $S^{2n}$, say, as
the upper hemisphere, and extend $f$ from the disc  to the whole
sphere by identity. It is rather straightforward to check that,
once orientations are fixed, the isotopy class of $\hat f$ does
not depend on the choices made. It is known  that if $\Sigma_f$ is
an exotic sphere, then the connected sum $\Cal T^{2n+1}_f = \Bbb
T^{2n+1}\#\Sigma_f$ is an exotic torus, i.e. it is topologically,
but {\it not smoothly} homeomorphic to $\Bbb T^{2n+1}$. Moreover,
it fibers over $S^1$ with the gluing map $f,$ see \cite{7},4.3.
This implies that $ \Cal T^{2n+2}=(\Bbb T^{2n+1}\#\Sigma_f)\times
S^1$ is an exotic torus (cf. Sec. 4) which fibers over $\Bbb T^2$
with fiber $\Bbb T^{2n}$.

Note that if $f$ were isotopic to a symplectomorphism, this fibration
would be symplectic. Clearly,  the cohomology ring
$H^*(\Cal T^{2n+2})$ is isomorphic to $H^*(\Bbb T^{2n})\otimes H^*(\Bbb T^2)$.
Hence, the cohomological  condition of Thurston's construction
were also satisfied and we would get a symplectic structure
on the exotic torus $\Cal T^{2n+2}_f$. These observations motivate the following questions.
\proclaim{Problem 1} Is there a symplectic structure on a torus with an
exotic differential structure?
\endproclaim
\proclaim{Problem 2} Is the fibration $\Bbb T^{2n}\to\Cal T^{2n+2}\to \Bbb T^2$
symplectic?
\endproclaim
Equivalently, one can ask

\proclaim{Problem 3} Given a diffeomorphism $f:\Bbb T^{2n}\to\Bbb T^{2n}$
supported in  an embedded disc but
non-isotopic to the identity, is there a symplectomorphism in the isotopy class of $f$?
\endproclaim
Our goal is to give negative examples to Problem 3  in the
particular case of the {\it standard symplectic structure on $\Bbb
T^{2n}$}. Note that this problem is related to the question posed
by McDuff and Salamon \cite{10}, p. 328, whether every
symplectomorphism of a torus which acts trivially on homology is
isotopic to the identity.

Let $\pi_0(\operatorname{Diff_+}\,(M))$ denote the group of
isotopy classes of orientation preserving  diffeomorphisms  of a
smooth oriented manifold $M$. Assume now that $M$ is
$2n$-dimensional and admits almost complex structures, and let
${\Bbb J}M$ denote the set of homotopy classes of such structures,
compatible with the given orientation. Any diffeomorphism $f$ acts
on the set of all almost complex structures by the rule
$$f_*J=df J df^{-1},$$
where $df: TM\to TM$ denotes the differential of $f$. This action
clearly descends to the action of
$\pi_0(\operatorname{Diff_+}\,(M))$ on ${\Bbb J}M$.

Let now ${\frak G}\, (M)$ denote the subgroup of $\pi_0(\operatorname{Diff}\,(M))$
generated by diffeomorphisms with supports in discs.

We have seen above that the positive answer to Problem 3 would
imply the positive answer to Problems 1 and 2. In the sequel we
will show that the answer to Problem 3 is {\it negative} in case
$2n\equiv 0$ (mod \,8), and the standard symplectic structure
$\omega_0$: there exist diffeomorphisms $f: \Bbb T^{8k}\to \Bbb
T^{8k}$ supported in a disc, whose isotopy classes $[f]\in\frak
G\,(M)$ do not preserve the homotopy class $[J_0]\in{\Bbb J}M$ of
the standard complex structure. Therefore,  they cannot be
isotopic to symplectomorphisms with respect to the standard
symplectic structure $\omega_0$. Indeed, given a symplectic form
$\omega ,$ any diffeomorphism preserving $\omega$ preserves also
the homotopy class of  any almost complex structure compatible
with $\omega$ (since  the space of all such almost complex
structures is contractible, hence connected).

Now, let us describe in more detail the results of the paper. We
will obtain first a necessary homotopic condition on a
diffeomorphism to be isotopic to a symplectomorphism.

\proclaim{Theorem 1} Let $f \in Diff\, (\Bbb T^{4n})$ be supported
in a disc $D^{4n}\subset \Bbb T^{4n}.$ If $f$  is isotopic to a
symplectomorphism with respect to the standard symplectic
structure, then $df$ restricted to its support disc $D^{4n}$ gives
in $\pi_{4n}SO(4n)$ the  trivial  homotopy class.
\endproclaim

Remark. The same proof is valid for any symplectic form compatible
with a parallelizable (as {\it complex} vector bundle structure on
the tangent bundle) almost complex structure. It is not known
whether $\Bbb T^{2n}$ admits a non-parallelizable symplectic
structure. For example, for any almost complex structure on $\Bbb
T^4$ compatible with a symplectic structure its first Chern class
vanishes and no examples of symplectic forms on $\Bbb T^4$
non-homotopic to the standard form are known.

For $\Cal T^{8k+1}_f=\Bbb T^{8k+1}\#\Sigma_{\hat f}$ consider the
Atiyah - Milnor - Singer  invariant, alias the $\hat a$-genus
$\hat a$  \cite{1,7} (see Section 4). Since the $\hat a$-genus is
additive with respect to the operation of connected sum and it is
nontrivial for some homotopy spheres in dimension $8k+1,$ one
easily concludes that there exist isotopy classes of
diffeomorphisms $f$ with support in a disc such that $\hat a(\Cal
T^{8k+1}_f)\ne 0$. On the other hand, we prove the following
result.

\proclaim{Theorem 2} In the notation of Theorem 1, if $[df]=0$,
then $\hat a(\Cal T_f)=0$.
\endproclaim
Comparing Theorem 1 and Theorem 2 we come to the main conclusion
of the paper.

\proclaim{Theorem 3} For any $k>0$ there exist diffeomorphisms
$f:\Bbb T^{8k}\to \Bbb T^{8k}$ with support in a disc which are
not isotopic to a symplectomorphism of $(\Bbb T^{8k}, \omega_0).$
\endproclaim
Of course, Theorem 3 implies all observations we have mentioned, and, in
particular, a partial negative answer to problem 2 in case of the standard
symplectic torus $\Bbb T^{2n}$.

The paper is organized as follows. Sections 2 and 3 are devoted to
the proof of Theorem 1, Section 4 describes the technique of the
$\hat a$-genus, and, finally, Section 5 contains the proofs of
Theorems 2 and 3. \vskip6pt {\bf Acknowledgement.} We are grateful
to Dusa McDuff for valuable advice and to Yuli Rudyak for reading
the manuscript and several useful remarks. We are also indebted to
the referee for pointing out a gap in our previous version of
Lemma 3.

\vskip6pt

\head 2. Preparatory work for the proof of Theorem 1
\endhead
To prove Theorem 1, we need some preparatory work.
\subhead 2.1 Interpretation of $[df]$
\endsubhead
Consider  diffeomorphisms
$f : D^m \to D^m $
 which are equal to the identity near the boundary.
The group of such diffeomorphisms we denote by
$\operatorname{Diff}_{\epsilon}(D^m,\partial D^m).$
The differential $df$ of such
a diffeomorphism yields a map
$(D^m,\partial D^m)\to (GL(m,\Bbb R),\{ id\} )$ when we use
a trivialization of the tangent bundle of the disc.
Since such a trivialization is unique up to
homotopy, we get a well defined  homotopy class
$[df]\in \pi_mGL(m,\Bbb R) \cong \pi_mSO(m)$.

Given an embedding of $D_m$ in a manifold, we can always
extend $f$ to the whole manifold,
since $f=\operatorname{id}$ near the boundary of the disc.
In particular,  $f$ extended by the identity map $\operatorname{id}$
to $S^m$ gives the  diffeomorphism
$$\hat f: S^m=D^m\cup_{\partial D^m}D^m\to S^m=D^m\cup_{\partial D^m}D^m.$$

If $M$ is an arbitrary manifold, and $f\in{\frak G}(M)$ is a
diffeomorphism with support in a disc, the restriction of $f$ to
this disc yields as above a homotopy class in $\pi_mSO(m)$. For
simplicity, we will denote this class by the same symbol
$[df]\in\pi_mSO(m)$. \subhead 2.2 Homomorphism $\Gamma_{m}\to
\pi_{m-1}SO(m-1)$
\endsubhead
Let us denote $\Gamma_{m}=\pi_0\operatorname{Diff_+}S^{m-1},$
where $\operatorname{Diff}_+S^k$ denotes the group of orientation
preserving diffeomorphisms of the $k$-sphere. It is well known
that $$\Gamma_m = \pi_0\operatorname{Diff_+}S^{m-1}\cong
\pi_0\operatorname{Diff}_{\epsilon}(D^{m-1},\partial D^{m-1}).$$
It follows that there is a well defined homomorphism
$$
\CD
\Gamma_m@>{\cong}>>\pi_0\operatorname{Diff}_{\epsilon}(D^{m-1},\partial
D^{m-1}) @>{\mu}>>\pi_{m-1}SO(m-1),
\endCD
$$
where
$$\mu([f])=[df].\eqno (1)$$
\subhead 2.3 Stabilization maps
\endsubhead
Consider the natural inclusions
$$SO(m-1)\hookrightarrow SO(m)\hookrightarrow SO(m+1)$$
and the induced maps of homotopy groups
$$
\CD
\pi_{m-1}SO(m-1) @>{j_m}>>\pi_{m-1}SO(m) @>{J_m}>> \pi_{m-1}SO(m+1)
\endCD
$$
which we will call the {\it stabilization maps}. From the work of
Kervaire \cite{8} one can obtain all the homotopy groups
$\pi_{m-1}SO(m-1),\pi_{m-1}SO(m),\pi_{m-1}SO(m+1)$ for all $m$, as
well as the kernels of the stabilization maps $J_m$. In the sequel
we will need this information. Because of that, we give  the table
of the homotopy groups together with $\operatorname{Ker}J_m$ up to
$m=15.$ It is well known that for $m>7$ the homotopy groups we are
interested in are 8-periodic. \vskip10pt

\centerline{\bf Table 1}

\vskip 0.5cm \vbox{\tabskip=0pt\offinterlineskip
\def\tablerule{\noalign{\hrule}}
\halign to 400pt{\strut#&\vrule#\tabskip=1em plus2em& \hfil# \hfil
& \vrule#&\hfil# \hfil & \vrule#&\hfil# \hfil & \vrule#&
\hfil#\hfil & \vrule#&\hfil { #} \hfil & \vrule # \tabskip=0pt\cr
\tablerule\omit&height5pt&\omit&height5pt&\omit&height5pt&\omit&height5pt
&\omit&height5pt&&\cr&&

 m && $\pi_{m-1}SO(m-1)$ && $\operatorname{Ker}J_m$ && $\pi_{m-1}SO(m)$
 && $\pi_{m-1}SO(m+1)$ &\cr
\tablerule\omit&height5pt&\omit&height5pt&\omit&height5pt&\omit&height5pt
&\omit&height5pt&&\cr && 3 && 0 && 0 && 0 && 0 &\cr
\tablerule\omit&height5pt&\omit&height5pt&\omit&height5pt&\omit&height5pt
&\omit&height5pt&&\cr && 4 && $\Bbb Z$ && $\Bbb Z$ && $\Bbb
Z\oplus\Bbb Z$ && $\Bbb Z$ &\cr
\tablerule\omit&height5pt&\omit&height5pt&\omit&height5pt&\omit&height5pt
&\omit&height5pt&&\cr && 5 && $\Bbb Z_2\oplus\Bbb Z_2$ && $\Bbb
Z_2$ && $\Bbb Z_2$ && 0 &\cr
\tablerule\omit&height5pt&\omit&height5pt&\omit&height5pt&\omit&height5pt
&\omit&height5pt&&\cr && 6 && $\Bbb Z_2$ && $\Bbb Z$ && $\Bbb Z$
&& 0 &\cr
\tablerule\omit&height5pt&\omit&height5pt&\omit&height5pt&\omit&height5pt
&\omit&height5pt&&\cr && 7 && 0 && 0 && 0 && 0 &\cr
\tablerule\omit&height5pt&\omit&height5pt&\omit&height5pt&\omit&height5pt
&\omit&height5pt&&\cr && 8 && $\Bbb Z$ && $\Bbb Z$ && $\Bbb
Z\oplus\Bbb Z$ && $\Bbb Z$ &\cr
\tablerule\omit&height5pt&\omit&height5pt&\omit&height5pt&\omit&height5pt
&\omit&height5pt&&\cr && 9 && $\Bbb Z_2\oplus\Bbb Z_2\oplus\Bbb
Z_2$ && $\Bbb Z_2$ && $\Bbb Z_2\oplus\Bbb Z_2$ && $\Bbb Z_2$ &\cr
\tablerule\omit&height5pt&\omit&height5pt&\omit&height5pt&\omit&height5pt
&\omit&height5pt&&\cr && 10 && $\Bbb Z_2\oplus\Bbb Z_2$ && $\Bbb
Z$ && $\Bbb Z\oplus\Bbb Z_2$ && $\Bbb Z_2$ &\cr
\tablerule\omit&height5pt&\omit&height5pt&\omit&height5pt&\omit&height5pt
&\omit&height5pt&&\cr && 11 && $\Bbb Z_4$ && $\Bbb Z_2$ && $\Bbb
Z_2$ && 0 &\cr
\tablerule\omit&height5pt&\omit&height5pt&\omit&height5pt&\omit&height5pt
&\omit&height5pt&&\cr && 12 && $\Bbb Z$ && $\Bbb Z$ && $\Bbb
Z\oplus\Bbb Z$ && $\Bbb Z$ &\cr
\tablerule\omit&height5pt&\omit&height5pt&\omit&height5pt&\omit&height5pt&
\omit&height5pt&&\cr && 13 && $\Bbb Z_2\oplus\Bbb Z_2$ && $\Bbb
Z_2$ && $\Bbb Z_2$ && 0 &\cr
\tablerule\omit&height5pt&\omit&height5pt&\omit&height5pt&\omit&height5pt&
\omit&height5pt&&\cr && 14 && $\Bbb Z_2$ && $\Bbb Z$ && $\Bbb Z$
&& 0 &\cr
\tablerule\omit&height5pt&\omit&height5pt&\omit&height5pt&\omit&height5pt&
\omit&height5pt&&\cr && 15 && $\Bbb Z_4$ && $\Bbb Z_2$ && $\Bbb
Z_2$ && 0 &\cr \tablerule\noalign{\smallskip}\hfil\cr}}

\subhead
 2.4 Image $j_m([df])$
\endsubhead
We will  need the following result.
\proclaim{Lemma 1} The following equality holds for any $m$:
$$j_m([df])=0.$$
\endproclaim
\demo{Proof} The proof is a consequence of the relations
between homotopy groups from  the Table. The necessary calculation falls into the cases:
\roster
\item "(i)" $m\equiv 4$\, (mod 8),
\item "(ii)" $m$ is odd, $m\ne 3,7$,
\item "(iii)" $m\equiv 2$\, (mod 8), or $m\equiv 6$\, (mod 8),
\item "(iv)" $m=3,7$,
\item "(v)" $m\leq 4.$
\endroster
Let us start with the following observation. Consider the homotopy
sphere $\Sigma_f^m=D^m\cup_fD^m$ and the tangent bundles $TS^m$
and $T\Sigma^m_f$. Since, up to homotopy, $m$-dimensional vector
bundles over homotopy spheres
are classified by elements of $\pi_{m-1}SO(m)$ %\cite{St},
one can check that, if $\tau\in\pi_{m-1}SO(m)$ represents $TS^m$,
then the element $\tau+j_m([df])$ represents $T\Sigma_f^m$.
Furthermore, the stable tangent bundle of $\Sigma_f^m$ is trivial.
For $m>4$ it is proved in \cite{9} and for $m\leq 4$ the group
$\Gamma_m$ is trivial by \cite{12}, \cite{5}. In particular,
$j_m([df])\in\operatorname{Ker}\,J_m.$

Again by Kervaire and Milnor \cite{9}, the group $\Gamma_m$ is
finite for $m>5,$ since it is isomorphic to the group of
h-cobordism classes of homotopy spheres (cf. Cerf \cite{6}). Now
we are ready to exhibit our case by case calculation. \roster
\item "Case (i)" For $m\equiv 4$\, (mod 8), the group
$\pi_{m-1}SO(m-1)$ is torsion free (see Table 1), while $\Gamma_m$
is always finite. Thus $[df]=\mu([f])$ must be zero. \item "Case
(ii)" Here the table yields $\operatorname{ker}\,J_m=\Bbb Z_2$.
Since the vector bundles $TS^m$ and $T\Sigma_f^m$ are both
nontrivial, we see that both $\tau$ and $\tau+j_m([df])$ are
non-zero elements of $\Bbb Z_2$. The only possibility for that is
$\tau=\tau+j_m([df])$, and again $j_m([df])=0.$ \item "Case (iii)"
Using Table 1 again, one can notice that $\pi_{m-1}SO(m-1)$ is
finite, while $\operatorname{Ker}\,J_m$ is torsion free. Hence,
$j_m([df])\in\operatorname{Ker}\,J_m$ can be only zero. \item
"Case (iv)" If $m=3,7$, then $\operatorname{Ker}\,J_m=0$, and
there is nothing to prove. \item "Case (v)" For $m\leq 4$ the
group $\Gamma_m$ is trivial.
\endroster
The proof is complete.

\hfill$\square$
\enddemo
\proclaim{Corollary 1} For any homotopy  sphere $\Sigma $ of
dimension $n$ its tangent bundle is isomorphic to the tangent
bundle of $S^n.$
\endproclaim

This implies for example the equality $span\, \Sigma = span\,
S^n$. Here $span$ denotes the maximal number of linearly
independent vector fields. This equality was proved (for any
stably parallelizable manifold) in \cite{4}, \cite{13} by
different arguments.

\subhead 2.5 $\pi_0(\operatorname{Diff}(M))$-action on ${\Bbb J}M$
\endsubhead
If $M^{2n}$ is a parallelizable smooth manifold  then any choice
of an almost complex structure $J_0$ on $M$ with a complex
trivialization of $(TM, J_0)$ determines a bijection
$$
\CD
{\Bbb J}M @>{\cong}>> [M, SO(2n)/U(n)]
\endCD
$$
such that $J_0$ corresponds to the class of constant map (cf.
\cite{11}, Prop.2.48). The correspondence can be briefly described
as follows. Any two complex structures on $\Bbb R^{2n}$ compatible
with a given orientation are equivalent up to a linear,
orientation preserving isomorphism. For a manifold, any local
complex trivialization of $(TM,J_0)$ give locally an
identification of complex structures on $T_xM$ with $GL_+(2n,\Bbb
R)/GL(n,\Bbb C).$ Globally, any almost complex structure
corresponds to a (smooth) section of the bundle with fiber
$GL_+(2n,\Bbb R)/GL(n,\Bbb C),$ associated to the tangent bundle.
If the manifold is parallelizable, then any fixed complex
parallelization provides a 1-1 correspondence between the set of
such sections  and the set of maps $M \>>> GL_+(2n,\Bbb
R)/GL(n,\Bbb C).$  Passing to homotopy classes and replacing the
target space by its homotopy equivalent $SO(2n)/U(n)$ we get $[M,
SO(2n)/U(n)].$ Note that taking the parallelization we assume
$J_0$ to give the trivial complex structure on the tangent bundle
of $M^{2n},$ but other structures can be arbitrary. The
correspondence above measures the difference between any almost
complex structure and the given one. In particular $J_0$ is
sent to the class of the constant map, and if $J$ is trivial, then
there exists a map $g:M\>>> GL_+(2n,\Bbb R)$ such that $J=gJ_0g^{-1}.$
For a nontrivial $J$ such a map exists only locally. However, different
local choices differ by a function to $GL(n,\Bbb C),$ thus
give a well-defined map to $GL_+(2n,\Bbb R)/GL(n,\Bbb C).$

Using such correspondence we readily see the action of
$Diff_+(M)$ on $\Bbb J(M)$ as follows. Given $f$ and $J,$
represent $[df]$ by a map $\delta : M\>>> SO(2n)$ and $[J]$ by
$\phi :M\>>> SO(2n)/U(n).$ Then we have
$[f]_*[J]=[\delta\centerdot \phi ]\in [M,SO(2n)/U(n)],$  where
$\delta\centerdot \phi$ is given by the natural action of $SO(2n)$
on $SO(2n)/U(n).$

Consider the bundle
$$U(n)\to SO(2n)\to SO(2n)/U(n)$$
and the corresponding part of the long homotopy exact sequence
$$... \to \pi_{2n}U(n)\to\pi_{2n}SO(2n)\to \pi_{2n}SO(2n)/U(n)\to...$$
Our first goal is to prove that if
$[f]\in\pi_0(\operatorname{Diff}(\Bbb T^{8k}))$ is determined by a
diffeomorphism with support in a disc, then it acts  on ${\Bbb
J}\Bbb T^{8k}$ with fixed points if and only if
$[df]\in\operatorname{Im}(\pi_{8k}U(4k))\subset\pi_{8k}SO(8k).$
The proof requires two technical facts we will give now.

\proclaim{Lemma 2} Let $D^{2n}\subset M^{2n}$ be an embedded disc
of codimension zero and let $\theta$ denote the map $M\to S^{2n}$
obtained by shrinking the complement
$M\setminus\operatorname{Int}\,D^{2n}$ to a point:
$$\theta: M\to S^{2n}=\operatorname{Int}\,D^{2n}
\cup\{pt\},\,\theta(M\setminus\operatorname{Int}
\,D^{2n})=\{pt\},\,\theta|_{\operatorname{Int}\,D^{2n}}=\operatorname{id}.$$
 If $M=\Bbb T^{2n}$, then the induced map of the sets of homotopy classes of maps
$$
\CD
[S^{2n},SO(2n)/U(n)] @>{\hat\theta}>> [M,SO(2n)/U(n)]
\endCD
$$
given by $[\psi]\to [\psi\circ\theta ],\,[\psi]\in [S^{2n},SO(2n)/U(n)]$, is injective.
\endproclaim
\demo{Proof} Let $M'$ denote the complement of
$\operatorname{Int}D^{2n}$ and let $X=SO(2n)/U(n)$ (in this
particular argument the target space $X$ may be arbitrary).
Consider  the Puppe long exact sequence of homotopy classes
(\cite{15}, 6.21)
$$... @>>> [SM, X] @>>> [SM',X] @>>> [M/M',X]@>>>[M,X] @>>>[M',X],$$

where $S$ denotes the unreduced suspension. Note that $X$ is
simply connected so we do not care about basepoints. It is clear,
that $[M/M',X]\to [M,X]$ can be identified with the map
$\hat\theta: [S^n,X]\to [M,X]$, and, therefore, the injectivity of
$\hat\theta$ will follow, if one proves the {\it surjectivity} of
the term
$$[SM, X]@>>> [SM',X]$$
in the Puppe long exact sequence.
The latter is naturally identified with the map
$[SM,X]\to [S(M\setminus\{pt\}),X]$ induced by the inclusion
$$i: S(M\setminus\{pt\})\hookrightarrow SM.$$
Now, take into consideration the homotopy equivalence
$$S\Bbb T^{2n} \simeq\vee S^j.\eqno (2)$$
The latter can be explained as follows. Use the well-known formulae
$$S(X\vee Y)=SX\vee SY\,\,\text{and}\,\,S(X\wedge Y)=X
\wedge Y\wedge S^1=SX\wedge Y=X\wedge SY$$
and
$$S^1\wedge X=SX,$$
where $X\wedge Y$ denotes the smash-product of spaces $X$ and $Y$.
Applying these formulae to $\Bbb T^{2n}=S^1\times ...\times S^1$
one obtains $(2)$.
 In the latter equivalence, the number of  spheres of
dimension $j$ equals to the Betti number $b_{j-1}$ of $\Bbb
T^{2n}$. We have the following (homotopy) commutative diagram
$$
\CD
S(\Bbb T^{2n}\setminus\{pt\}) @>>> \vee S^j\setminus S^{2n+1}\\
@V{i}VV @V{i_S}VV\\
S\Bbb T^{2n} @>>> \vee S^j
\endCD
$$
where $i_S$ is an obvious map of the wedges of spheres. Indeed, if
one suspends the CW-complex $\Bbb T^{2n}$ with a standard cell
decomposition, the suspension of every $k$-cell will give a
$(k+1)$-dimensional sphere, and therefore, cutting the top cell
and suspending the rest will give the wedge of spheres with the
top one shrunken to a point. It follows that there is a retraction
$r: S\Bbb T^{2n}\to S(\Bbb T^{n}\setminus\{pt\})$ corresponding to
the right vertical arrow of the previous diagram, and given by
shrinking the sphere of the maximal dimension to the point. It
follows that the induced map $i_*=(i_S)_*: [S\Bbb T^{2n}, X]\to
[S(\Bbb T^{2n}\setminus\{pt\}),X]$ has a right inverse $r_*$. This
is the same as the surjectivity of $i_*$, as required.
\hfill$\square$
\enddemo

\proclaim{Lemma 3} Consider a diffeomorphism $f$ of  $\Bbb T^{2n}$
with support in a disc and a parallelizable almost complex
structure $J_0$ on $\Bbb T^{2n}.$ Then $f$ preserves the homotopy
class of $J_0$ if and only if
$[df]\in\operatorname{Im}(\pi_{2n}U(n))\subset \pi_{2n}SO(2n)$.
\endproclaim
\demo{Proof} Let $*$ denotes the base point of $SO(2n)/U(n)$
corresponding to $U(n)$ and $D_0$ the disc in $\Bbb T^{2n}$
containing the support of $f.$ By {\it the constant map} we will
understand here a map which sends every point to $*.$ As we have
explained above, the correspondence $\Bbb J(\Bbb T^{2n}) \cong
[\Bbb T^{2n}, SO(2n)/U(n)]$ can be chosen such that $J_0$
corresponds to the class of the constant map, denoted $\phi .$ Let
$df$ be represented by $\delta : \Bbb T^{2n}\rightarrow SO(2n),$
constant outside $D_0,$ and assume $[f^*J_0]=[J_0].$ This gives a
homotopy $H:\delta = \delta\centerdot \phi \sim \phi.$ We have to
prove that there is another homotopy which is constant outside
$D_0.$ The homotopy $H$ gives a map $h:\Bbb T^{2n}-D_0\>>>
\Omega_*SO(2n)/U(n)$ to the space of loops at $*.$ Every subtorus
$\Bbb T^k\subset \Bbb T^{2n}$ is the image of a retraction $r:\Bbb
T^{2n}\>>> \Bbb T^k,$ thus any map on $\Bbb T^k$ extends to $\Bbb
T^{2n}.$ In particular we have a map $\Phi =hr:\Bbb T^{2n}\>>>
\Omega_*SO(2n)/U(n)$ equal to $h$ on $\Bbb T^k.$ Inverting $\Phi$
pointwisely and composing the given homotopy $H$ with  the
resulting homotopy $\Phi^{-1}: \phi \sim \phi$ we get a new
homotopy $\delta\sim \phi$ which is homotopy trivial on $\Bbb
T^k.$ Thus there is also a homotopy $\hat f \sim \phi$ which is
constant on a neighborhood of $T^k.$

Let $\Bbb T^{2n} = S^1\times ...\times S^1,$ and denote by $N_i$ a
tubular neighborhood of $\Bbb T_i = S^1\times ...\times \{ x_i\}
\times ...\times S^1,$ where $x_i$ is a point in $i-$th $S^1$
factor. Then $\Bbb T^{2n}-\bigcup_{i=1}^{2n}N_i$ is a disc and we
can assume that it is equal to $D_0.$  Applying inductively  the
above construction to all subtori $\Bbb T_i$ we get eventually  a
homotopy $\delta \sim \phi$ which is constant in the complement of
$D_0.$ Thus the image of $[df]$ in $\pi_{2n}SO(2n)/U(n)$ is zero,
and the exactness of the homotopy sequence of the fibration
$U(n)\rightarrow SO(2n)\rightarrow SO(2n)/U(n)$ yields our claim.

 \hfill$\square$
\enddemo

Remark. It is conceivable that Lemma 3 is valid without assumption
that $J_0$ is parallelizable.

\head 3. Proof of Theorem 1
\endhead
\noindent By Lemma 3, if $[f]$ is isotopic to a symplectomorphism,
then $[df] \in \operatorname{Im}(\pi_{4n}U(2n))$. We will show,
that if $[df]$ is nonzero in $\pi_{4n}SO(4n)$, then $[df]\notin
\operatorname{Im}(\pi_{4n}U(2n))$, and this will complete the
proof.

The argument goes as follows. Both $U(2n+1)$ and $SO(4n+2)$ act
transitively on $S^{4n+1}$ and yield  diffeomorphisms
$$S^{4n+1}\cong SO(4n+2)/SO(4n+1)\cong U(2n+1)/U(2n).$$
Consider then the following commutative diagram of fibrations and natural inclusions
$$\CD
SO(4n+1) @>>> SO(4n+2) @>>> S^{4n+1}\\
@A{j\circ i}AA @AAA @A{\operatorname{id}}AA \\
U(2n)  @>>> U(2n+1) @>>> S^{4n+1},
\endCD
$$
where $i: U(4n)\hookrightarrow SO(4n)$ and $j:
SO(4n)\hookrightarrow SO(4n+1)$ denote the natural inclusions.
This diagram
 yields the commutative diagram of the group homomorphisms
$$
\CD
\pi_{4n+1}S^{4n+1} @>{\partial}>> \pi_{4n}SO(4n+1)@>{J_{4n+1}}>>\pi_{4n}SO(4n+2)\\
@A{\operatorname{id}}AA @A{j_{4n+1}\circ i_*}AA @AAA\\
\pi_{4n+1}S^{4n+1} @>{\partial '}>> \pi_{4n}U(2n)
@>>>\pi_{4n}U(2n+1)
\endCD
$$
where the horizontal rows represent parts of the long homotopy
sequences of the corresponding fibrations and $\partial, \partial
'$ denote the connecting homomorphisms.
%?? niepotrzebne?Here $i_*$ denotes the map
%of homotopy groups induced by $i$.
Let $\alpha$ and $\beta$
denote, respectively, the  generators of cyclic groups (see
\cite{3}):
$$\pi_{4n+1}S^{4n+1}\cong\Bbb Z\cong\langle\beta
\rangle,\,\,\text{and}\,\,\pi_{4n}U(2n)\cong\Bbb
Z_{(2n)!}\cong\langle\alpha\rangle.$$
By \cite{8, Lemma I.1}
$$\alpha=\partial '\beta.$$
Now, the kernel $\operatorname{Ker}\, J_{4n+1}$ is nontrivial and
$$\operatorname{Ker}\,J_{4n+1}=\Bbb Z_2.$$
It follows that
$$\partial (\pi_{4n+1}S^{4n+1})=\operatorname{Ker}\,J_{4n+1}=\Bbb Z_2.$$
In particular, $\partial\beta\ne 0$. From the commutativity of the diagram
$$0\ne \partial\beta=(j_{4n+1}\circ i_*)\partial '\beta=j_{4n+1}(i_*\alpha).$$
It follows that $i_*\alpha\notin\operatorname{Ker}\,j_{4n+1}.$ In
fact, it gives more: since $\pi_{4n}SO(4n)$ is a 2-torsion group
and $\pi_{4n}U(2n)$ is cyclic, it follows that
$\operatorname{Im}\,i_*\cap\operatorname{Ker}\,j_{4n+1}=\{ 0\}.$
But from Lemma 1 we know that
$[df]\in\operatorname{Ker}\,j_{4n+1},$ thus
$[df]\in\operatorname{Im}(\pi_{4n}U(2n))$ only if it is zero, as
required.

\hfill$\square$

\head 4. $\hat a$-genus and family index
\endhead
The  tool that yields an obstruction to isotoping the
diffeomorphism $f$ (supported in a disc) to a symplectomorphism is
a $KO$-theoretical invariant, namely, the $\hat a$-genus of a
closed spin manifold. In this section we summarize the properties
of the $\hat a$-genus and explain  the argument. A beautiful
presentation of  techniques used in this section one can find in
the monograph of Lawson and Michelson \cite{11}.

 We consider closed spin manifolds. It is known that  stably parallelizable
manifolds are spin, hence
$\Bbb T^m$ and $\Sigma_f$ are. Also, connected sums of spin
manifolds admit spin-structures, which implies that
${\Cal T}_f=\Bbb T^m\#\Sigma_f$ is spin.

The $\hat a$-genus can be defined as the $KO$-theoretical index of
the Dirac operator. It is known that the coefficient groups
$KO^{-*}(pt)$ are the following
$$
KO^{-m}(pt)= \cases
\Bbb Z &\,\,\text{if}\,\,m \equiv 0 \, \, (\text{mod}\, 8);\\
\Bbb Z_2 &\, \, \text{if}\, \, m\equiv 1,2\, \,  (\text{mod}\, 8);\\
0 &\,\,\, \text{for any other} \,\,\,m.
\endcases
$$

Let $f: M\to \{pt\}$ denote the obvious collapsing map.
For a spin structure on $M$ we have the Gysin (or direct image) map
$$f_{!}: KO^0(M)\to  KO^{-m}(pt),$$
where $m=\dim\, M.$

\definition{Definition} The $\hat a$-genus of a closed
spin manifold $M^m$ is an element of \newline $KO^{-m}(pt)$ given
by the formula
$$\hat a(M)=f_{!}(1).$$
\enddefinition

This is the topological index, with values in $KO^{-*}(pt),$ of
the Dirac operator on $M.$ In fact we use only the torsion part of
$\hat a$ called  the {\it Hitchin invariant}.

Now we formulate explicitly the properties of $\hat a$ used in the
sequel.

\proclaim{Proposition 1} The $\hat a$-genus has the following
properties: \roster \item "(i)" for any closed spin manifolds $X$
and $Y$
$$\hat a(X\#Y)=\hat a(X)+\hat a(Y), \,\,\text{and}\,\,
\hat a(X\times Y)=\hat a(X)\cdot \hat a(Y),$$
\item "(ii)" $\hat a$ is a spin cobordism invariant,
\item "(iii)" for any $m>2$ and any spin structure on  the standard torus $\Bbb T^m$
we have $\hat a(\Bbb T^m)=0.$
\endroster
\endproclaim
\demo{Proof} See \cite{11}. Note that $(iii)$ follows from $(i)$
and $(ii)$ since any spin structure on $\Bbb T^m$ is given as
product of spin structures on circles. Recall that $S^1$ has two
spin structures, one of them has $\hat a\not=0$, and does not
bound. However, the third power of the non-zero element of
$KO^{-1}(pt)$ vanishes in the ring $KO^{-*}(pt)$.
\enddemo

\hfill$\square$

\proclaim{Corollary 2} The manifold $\Cal T_f = (\Bbb
T^{8k+1}\#\Sigma_{\hat f})\times S^1$ is homeomorphic, but not
diffeomorphic to the standard torus $\Bbb T^{8k+2}$, if $\hat
a(\Sigma_f)\neq 0.$ In fact, the $\hat a$-genus of $\Cal T_f$ does
depend on the choice of the spin structure.
\endproclaim
\demo{Proof} Apply Proposition 1 and use the spin structure on
$S^1$ with non-vanishing $\hat a$.
\enddemo

\hfill$\square$

\remark{Remark} A classification of smooth structures on
topological tori is described in the chapter "Fake Tori" in
\cite{14}.
\endremark
Let there be given a smooth fiber bundle
$$
\CD
F @>>> M @>{\pi}>> B
\endCD
\eqno (3)
$$
with fiber and base being  spin manifolds. For any continuous
family of elliptic differential operators $P$ on fibers of such
bundle there is a well-defined {\it family index}

$$\operatorname{Ind}_mP  \in KO^{-m}(B),\ \ m=\dim F.$$
(see Atiyah and Singer \cite{1}). Consider the particular case of
the parametric Dirac operator $D$. We assume that the spin
structure on $M$ is  the one induced by spin structures on $B$ and
$F.$ Then we have the formula  \cite{1,7}.
$$\operatorname{Ind}_mD=\pi_{!}(1)\in KO^{-m}(B).$$

>From the functoriality  of the Gysin map we have

 \proclaim{Lemma 4} If the family index of the Dirac operator on the fiber bundle
(3) vanishes, then $\hat a(M)=0$.
\endproclaim

\hfill$\square$

Now we are going to describe a  condition ensuring  the vanishing
of the family index $\operatorname{Ind}_m\,D$ of  the parametric
Dirac operator. This is certainly known to experts, but we haven't
found any appropriate reference in the literature.

For a parallelizable manifold, a given trivialization of
the tangent bundle induces a spin structure and a Riemannian
metric on (the tangent bundle of) the manifold. We consider
a fiber bundle with a continuous family of  parallelizations
of fibers and corresponding spin structures and Dirac operators on fibers.

\proclaim{Lemma 5} For any closed spin manifold $F$ of dimension
$8k$ and any fiber  bundle $F\to M\to B$  which
 admits a fiberwise parallelization, the family index
 $\operatorname{Ind}_{8k}D$ vanishes.
\endproclaim
\demo{Proof} The topological index of the family is given as an
element of the Real K-theory of $B,$ $KR^{-8k}(B)\cong
KO^{-8k}(B).$ Given a parallelization of $F,$ the $KR$-symbol
class of the Dirac operator of $F$ is identified with an element
of $KR(F\times \Bbb R^k \times \Bbb R^k).$ This element is given
by a map of product bundles $F\times \Bbb R^k \times \Bbb R^k
\times S_+ \to F\times \Bbb R^k \times \Bbb R^k \times S_-$ where
$S_+, S_-$ are half-spin representations of the group $Spin(8k)$
and over a point $(x,u,v)$ the map is the multiplication by
$u+iv.$ In dimension $8k$ parallelizable manifolds have trivial
$\hat a$-genera, thus the symbol class of the Dirac operator on
$F$ becomes zero after passing to $KR^{-8k}(pt).$ In fact, any
parallelization yields a trivialization of the resulting Real
bundle. If we consider a fiberwise parallelization, then we obtain
the product of the above by $B.$ In particular, $Ind_{8k}D \in
KR^{-8k}(B)$ is zero.
\enddemo

\hfill$\square$

\head 5. Proofs of Theorems 2 and 3
\endhead

Lemmas 4 and 5 yield the following observation we can use to
complete the proofs of Theorems 2 and 3.

\proclaim{Proposition 2} Assume that  $M$ is fibred over a closed
spin manifold $B$ with closed spin fiber $F$ of dimension $8k.$ If
the fibration admits a fiberwise parallelization and the spin
structure on $M$ is the one induced from the parallelization and
the spin structure of $B,$ then  $\hat a(M)=0.$
\endproclaim

\hfill$\square$

\subhead Proof of Theorem 2
\endsubhead
Consider ${\Cal T}_f=\Bbb T^{8k+1}\#\Sigma_f$. We have already
mentioned that ${\Cal T}_f$ fibers over $S^1$ with $\Bbb T^{8k}$
as a fiber and the gluing map $\hat f.$

If $[df]=0,$ then the fibration admits fiberwise parallelization.
 Proposition 2 implies that $\hat a({\Cal T}_f)=0.$

\hfill$\square$

\subhead Proof of Theorem 3
\endsubhead
If $[f]$  preserves a homotopy class $[J]$ of an almost complex
structure, then $[df]=0$ and, by Theorems 1 and 2, $a({\Cal
T}_f)=0$. We have $\hat a({\Cal T}_f)=\hat a(\Bbb
T^{2n+1}\#\Sigma_f)=\hat a(\Sigma_f)$. However, it is well known
\cite{7} that in dimensions $8n+1$ there are homotopy spheres
$\Sigma_f$ with $\hat a(\Sigma_f)\not=0$. This completes the
proof.

\hfill$\square$

\remark{Remark} By Section 2, $\hat a({\Cal T}_f)$ is necessarily
a torsion element, so $\hat a$-genus can detect nontriviality of
$[df]$ only in dimensions $8n.$
However, the group where the class $[df]$ can take values is equal
 to $\Bbb Z_2$ for any even dimension.
We do not know  whether in  even dimensions $\neq 8n$
there are diffeomorphisms of spheres having
non-trivial homotopy class of the differential (while the
results of Section 2 show that $[df]=0$ if the dimension is odd).
\endremark

\medskip
\noindent Bogus\l aw Hajduk: {\bf Mathematical Institute, Wroc\l
aw University,

\noindent pl. Grunwaldzki 2/4, 50--384 Wroc\l aw, Poland}

\noindent \tt hajduk{\@}math.uni.wroc.pl

\noindent http://www.math.uni.wroc.pl/$\sim$hajduk

\medskip

\noindent Aleksy Tralle: {\bf Department of Mathematics and
Information

\noindent Tech\-no\-lo\-gy, University of Warmia and Mazury,}

\noindent Zolnierska 14A, 10-561 Olsztyn, Poland

\noindent \tt tralle{\@}matman.uwm.edu.pl

\noindent http://matman.uwm.edu.pl/$\sim$tralle

\vfill \eject


\Refs\nofrills{References} \widestnumber\key{1234}

\ref\key{1} \by M.F. Atiyah, I.M. Singer \paper The index of
elliptic operators IV,V \jour Ann. of Math. \vol 93 \yr 1971
\pages 119-138, 139-149
\endref\vskip6pt

\ref\key{2} \by C. Benson and C. Gordon \paper K\"ahler and
symplectic structures on nilmanifolds \jour Topology \vol 27 \yr
1988 \pages 513-518
\endref\vskip6pt

\ref\key{3} \by R. Bott \paper The stable homotopy of the
classical groups \jour Proc. Nat. Acad. Sci. USA \vol 43 \yr 1957
\pages 933-935
\endref\vskip6pt

\ref\key{4} \by G.E Bredon, A. Kosinski \paper Vector fields on
$\pi$-manifolds \jour Ann. of Math. \vol 84 \yr 1966 \pages 85 --
90
\endref\vskip6pt

\ref\key{5} \by J. Cerf \paper La nullit\'e de $\pi_0(DiffS^3)$
\jour Sem. Henri Cartan  1962/63 \pages Exp. 10-12, 20-21, Paris
1964
\endref\vskip6pt



\ref\key{6} \by J. Cerf \paper La stratification naturelle des
espaces de fonctions diff\'erentiables r\'eeles et le th\'eoreme
de la pseudoisotopie \jour Publ. IHES \vol 39 \yr 1970 \pages 5
--173
\endref\vskip6pt

\ref\key{7} \by N. Hitchin \paper Harmonic spinors \jour Adv.
Math. \vol 14 \yr 1974 \pages 1-55
\endref\vskip6pt

\ref\key{8} \by M. Kervaire \paper Some nonstable homotopy groups
of Lie groups \jour Illinois J. Math. \vol 4 \yr 1960 \pages
161-169
\endref\vskip6pt



\ref\key{9} \by M. Kervaire, J. Milnor \paper Groups of homotopy
spheres \jour Ann. of Math. \vol 77 \yr 1963 \pages 504-537
\endref\vskip6pt



\ref\key{10} \by H.B.Lawson, M.-L. Michelsohn \book Spin Geometry
\publ Princeton Univ. Press \yr 1989
\endref\vskip6pt

\ref\key{11} \by D. McDuff, D.Salamon \book Introduction to
Symplectic Topology \publ Oxford Univ. Press \yr 1998
\endref\vskip6pt




\ref\key{12} \by S. Smale \paper Diffeomorphisms of the 2-sphere
\jour Proc. of AMS \vol 10 \yr 1959 \pages 621 -- 626
\endref\vskip6pt

\ref\key{13} \by E. Thomas \paper Cross - sections of stably
equivalent vector bundles \jour Quart. J. of Math. \vol 17 \yr
1966 \pages 53 -- 57
\endref\vskip6pt

\ref\key{14} \by C.T.C. Wall \book Surgery on Compact Manifolds
(Second Edition, Edited by A.A. Ranicki) \publ American
Mathematical Society \yr 1999 \pages
\endref\vskip6pt

\ref\key{15} \by G.W. Whitehead \book Elements of Homotopy Theory
\publ Graduate Texts in Math. 61, Springer, Berlin \yr 1978
\endref\vskip6pt
\endRefs
\enddocument
\end